\documentclass[preprint,12pt]{elsarticle}
\usepackage{amssymb}
\usepackage{amsmath}

\newdefinition{rmk}{Remark}
\newproof{pf}{Proof}
\newproof{pot}{Proof of Theorem \ref{thm2}}
\author{Ting-Yang Hsiao \\ email: r04221011@ntu.edu.tw}
\journal{arXiv}

\begin{document}
\begin{frontmatter}

\title{The "Ant on a Rubber Rope" Paradox} 
\begin{abstract}
We clarify and generalize the \textit{ant on a rubber rope} paradox, which is a mathematical puzzle with a solution that appears counterintuitive. In this paper, we show that the ant can still reach the end of the rope even if we consider the step length of the ant and stretching length of the rubber rope as random variables.

\begin{keyword}
{Ant on a rubber rope; Law of large numbers; Paradox; Harmonic series} 
\end{keyword}
\end{abstract}
\end{frontmatter}

\section{Introduction}
The \textit{ant on a rubber rope} paradox is usually expressed as (see, for instance, [1] or [2]): \textit{An ant is at the left endpoint of a rubber rope, which is $1$ kilometre long and the ant crawls along the rope at a steady pace of $1$ centimetre per second. After the first second, the rubber rope stretches uniformly to $2$ kilometres instantly. After the next second, it stretches to $3$ kilometres, and so on. The question is, will the ant ever reach the right endpoint of the rubber rope?} 

Surprisingly, the answer to this paradox is "Yes". In the beginning, we should know that the ant moves with the rope when the rope is stretching. Now, we can analyse this paradox by measuring the fraction of the rope the ant covers after each second. When the sum of these fractions is $1$, the ant has come to the right endpoint of the rope.

Initially, the rubber rope is $1$ kilometre long and the ant moves $1$ centimetre forward. This means that at the end of the first second, the ant has travelled $\dfrac{1}{100,000}$ times of the rubber rope's length of $1$ kilometre. Meanwhile, after the first second, the rubber rope stretches uniformly to $2$ kilometres immediately. Then the ant crawls another centimetre after the next second. That is to say, the ant has gone an additional $\dfrac{1}{200,000}$ times of the rope's new length of $2$ kilometres. By repeating this process $m$ times, after the $m$th second, the ant's progress can be expressed as a fraction of the entire rubber rope, which is 
\begin{equation}\label{eqn: l0}
\dfrac{1}{100,000}\left( \dfrac{1}{1}+\dfrac{1}{2}+...+\dfrac{1}{m}\right).
\end{equation} 

The equation \eqref{eqn: l0} is a partial sum of the Harmonic series, which can be made as large as we desire. Therefore, the ant will complete the journey to the right endpoint of the rubber rope. 

Now, we consider the case that both the movements of the ant and the increments of the rubber rope are random. Specifically, let $X_0,X_1$,... be positive, independent, and identically distributed random variables with $EX_i=\mu_X>0$. Let $L_1,L_2$,... be positive, independent, and identically distributed random variables with $EL_i=\mu_L< \infty$. Let $l_0$ be a positive constant. From the very beginning, the length of the rope is $l_0$ units. At the first second, the ant, which is at the left endpoint of the rope, starts to move $X_0$ units. At the end of the first second, the rope uniformly stretches $L_1$ units to $l_0+L_1$ units immediately. In the next second, the ant moves another $X_1$ units. Again, at the end of the second second, the rope uniformly stretches another $L_2$ units to $l_0+L_1+L_2$ units. If this continues, will the ant still reach the other end of the rope? Astonishingly, the answer is still "Yes".\\ 

$Theorem~1$
In the stochastic model of the \textit{ant on a rubber rope}, the ant can still reach the other end of the rope almost surely.

In Section 2, we will establish two lemmas and then use them to prove Theorem 1. \\

\section{Two Lemmas}
To prove Theorem 1, we need the following Lemmas 1 and 2. 
\\
$Lemma~1$
An ant is originally at the left endpoint of a rubber rope of length $l_0$ units. At the $i$th second the ant moves along the rope at a pace of $x_{i-1}$ units, and then the rubber rope uniformly stretches $l_i$ units instantly. If, at the $m$th second, the ant still does not reach the right endpoint of the rope, then  the ant's progress can be expressed as a fraction of the entire rubber rope, which is 
\begin{equation}\label{eqn: 200}
\dfrac{x_0}{l_0}+\dfrac{x_1}{l_0+l_1}+\dfrac{x_2}{l_0+l_1+l_2}+...+\dfrac{x_{m-1}}{l_0+l_1+...+l_{m-1}}.
\end{equation} 
~\\

This lemma is not difficult to prove by induction argument. The key point is to note that the ratio of length of the ant's progress to that of the rope remains unchanged after the rope uniformly stretching.

$Lemma~2$
Let $X_0,X_1$,... be positive, independent, and identically distributed random variables with $EX_i=\mu_X>0$. Let $L_1,L_2$,... be positive, independent, and identically distributed random variables with $EL_i=\mu_L< \infty$. Let $l_0$ be a positive constant. Then \\
\begin{equation}
\dfrac{X_0}{l_0}+\dfrac{X_1}{l_0+L_1}+\dfrac{X_2}{l_0+L_1+L_2}+...=\infty
\end{equation}
almost surely. 


\begin{pf}[Proof]
Let $\Omega$ be the sample space. By strong law of large numbers there exists $\Omega_0$ with $P(\Omega_0)=1$, such that for each $\omega \in \Omega_0$, 
\begin{equation}
\dfrac{L_1(\omega)+...+L_n(\omega)}{n}\longrightarrow\mu_L,
\end{equation}
and
\begin{equation}
\dfrac{X_1(\omega)+...+X_n(\omega)}{n}\longrightarrow\mu_X,
\end{equation} 
when $n$ goes to $\infty$. This implies that for each $\omega \in \Omega_0$, and any given $\epsilon>0$, we can choose an $N \in \mathbb{N}$ such that
\begin{equation}\label{eqn: l3}
\left|\dfrac{l_0+L_1(\omega)+...+L_n(\omega)}{n}-\mu_L\right|<\epsilon,
\end{equation}
and
\begin{equation}\label{eqn: l1}
\left|\dfrac{X_1(\omega)+...+X_n(\omega)}{n}-\mu_X\right|<\epsilon,
\end{equation} 
if $n \geqslant N$. To simplify our notation, denote $L_i(\omega)$ and $X_i(\omega)$ by $l_i$ and $x_i$, respectively. \\
Now, for each $i\in \mathbb{N}$,
we have, according to \eqref{eqn: l1},
\begin{align}\label{eqn: l2}
\notag
&~~~\left|\dfrac{x_{(i-1)N+1}+...+x_{iN}}{N}-\mu_X\right|=\left|\dfrac{x_{1}+...+x_{iN}}{N}-\dfrac{x_{1}+...+x_{(i-1)N}}{N}-\mu_X\right| \\ \notag
&=\left|i\dfrac{x_{1}+...+x_{iN}}{iN}-i\mu_X+(i-1)\mu_X-(i-1)\dfrac{x_{1}+...+x_{(i-1)N}}{(i-1)N}\right| \\ 
\notag
&\leqslant i \left|\dfrac{x_{1}+...+x_{iN}}{iN}-\mu_X \right|+(i-1)\left|\dfrac{x_{1}+...+x_{(i-1)N}}{(i-1)N}-\mu_X\right| \\
&< i\epsilon+(i-1)\epsilon=(2i-1)\epsilon.
\end{align} \\
Consequently, for any fixed $m \in \mathbb{N}$, we can use \eqref{eqn: l3} and \eqref{eqn: l2} to yield that
{\renewcommand\baselinestretch{3.3}\selectfont
\begin{align}
\notag
&~~~\dfrac{x_0}{l_0}+\dfrac{x_1}{l_0+l_1}+\dfrac{x_2}{l_0+l_1+l_2}+...\\
\notag
&\geqslant \sum\limits_{i=1}^{N}\dfrac{x_i}{l_0+...+l_i}+\sum\limits_{i=N+1}^{2N}\dfrac{x_i}{l_0+...+l_i}+...+\sum\limits_{i=(m-1)N+1}^{mN}\dfrac{x_i}{l_0+...+l_i}\\
\notag
&\geqslant\sum\limits_{i=1}^{N}\dfrac{x_i}{l_0+...+l_N}+\sum\limits_{i=N+1}^{2N}\dfrac{x_i}{l_0+...+l_{2N}}+...+\sum\limits_{i=(m-1)N+1}^{mN}\dfrac{x_i}{l_0+...+l_{mN}} \\
\notag
&=\dfrac{\dfrac{x_1+...+x_N}{N}}{\dfrac{l_0+l_1+...+l_N}{N}}+\dfrac{1}{2}\dfrac{\dfrac{x_{N+1}+...+x_{2N}}{N}}{\dfrac{l_0+l_1+...+l_{2N}}{2N}}+...+\dfrac{1}{m}\dfrac{\dfrac{x_{(m-1)N+1}+...+x_{mN}}{N}}{\dfrac{l_0+l_1+...+l_{mN}}{mN}} \\
\notag
&\geqslant \dfrac{\mu_X-\epsilon}{\mu_L+\epsilon}+\dfrac{1}{2}\dfrac{\mu_X-3\epsilon}{\mu_L+\epsilon}+...+\dfrac{1}{m}\dfrac{\mu_X-(2m-1)\epsilon}{\mu_L+\epsilon} \rightarrow \dfrac{\mu_X}{\mu_L}\left( 1+\dfrac{1}{2}+...+\dfrac{1}{m}\right)~~\mbox{as}~\epsilon \rightarrow 0^+. 
\end{align}\par} 
Because the harmonic series diverges to $\infty$, we conclude that 
\begin{equation*}
\dfrac{X_0}{l_0}+\dfrac{X_1}{l_0+L_1}+\dfrac{X_2}{l_0+L_1+L_2}+...=\infty
\end{equation*}
almost surely. 
\end{pf} 

Now, we can prove Theorem 1 as follows. 

\begin{pf}[Proof of Theorem 1]
Using the notations as in Lemma 2, we define $$T=\min \left\{n ~\Bigg| ~\dfrac{X_0}{l_0}+\dfrac{X_1}{l_0+L_1}+...+\dfrac{X_{n-1}}{l_0+L_1+L_2+...+L_{n-1}}\geqslant 1 \right\}.$$
In view of Lemma 1, $T$ is just the time (in seconds) when the ant reaches the right endpoint of the rubber rope. Here we use the convention that $\min\emptyset=\infty$. Now, Lemma 2 tells that $T<\infty$ almost surely, which completes the proof.
\end{pf}




%
%
%
%

\end{document}